\documentclass[10pt,a4paper,english]{article}

\usepackage{amsfonts}
\usepackage{amssymb}
\usepackage{amsmath}
\usepackage{latexsym}
\usepackage{mathrsfs}

\usepackage{xcolor}	
	\usepackage{soul}

\usepackage{tikz}

\def\IN{\mathbb{N}}
\def\IZ{\mathbb{Z}}
\def\IQ{\mathbb{Q}}

\def\Otto{\Leftrightarrow}

\def\eps{\emptyset}
\def\sbsq{\subseteq}

\def\<{\langle}
\def\>{\rangle}

\def\[{\lbrack\!\lbrack}
\def\]{\rbrack\!\rbrack}

\newtheorem{lemma}{Lemma}[section]

\newtheorem{proposition}[lemma]{Proposition}

\title{The basis of Boole's logical theory}

\author{{Giovanna Corsi, Guido Gherardi}
\date{}
\\
\\
\small{Dipartimento di Filosofia e Comunicazione - Universit\`{a} di Bologna}
\\ \small{via Zamboni, 38 - I-40126 Bologna} 
\\ \small{giovanna.corsi@unibo.it} \quad  \small{guido.gherardi@unibo.it}  }

\begin{document} 

\maketitle

\begin{abstract}
 {In the present paper we aim to provide a thoughtful and exegetical account of the fundamental ideas at the basis of Boole's theory, with the goal of developing our investigation strictly within the conceptual structure originally introduced by Boole himself. In particular, we will focus on the meaning and the usefulness of the methods of the developments. We will also consider a slight variation of it that will allow us to to present in a new light some important and ingenuous aspects of Boole's calculus examined by the author in his work. Finally, a large attention is devoted to the analysis of the ``neglected'' logical connective of division.}
\end{abstract}

\section{\emph{Introduction and outline}}
Anyone approaching the study of George Boole's \emph{The Laws of Thought}, \cite{B1854}, a classic from the origins of modern logic, would expect to find in the critical literature a clear explanation of the fundamental ideas on which Boole based his theory of logic.  In our opinion no  such explanation has yet been given, though there are attempts to account for Boole's theory by means of complex algebraic systems, see \cite{H76} and \cite{BR09}.

 In the present paper we aim to provide a clear and exhaustive account of the fundamental ideas at the basis of Boole's whole theory of logic without appealing to anything other than concepts introduced by Boole himself {in the chapters I-VI of \emph{The Laws of Thought}.}



According to James W. van Evra \cite{VE77},
 ``While there is general agreement that his work occupies an important place in the history of logic, the exact nature of that importance remains elusive.''
We believe that we can show that the nature of that importance consists in ingenious  and refined ideas on which his logical theory is based.
The calculus that Boole takes as model for his logical theory is quantitative algebra. Both calculi are built on the same language based on the operations  $\times, +, - , / $; Boole goes to great lengths to show that those operations share the same formal properties since the variables of logical theory are interpreted on sets. Nevertheless, according to Boole, the formal correspondences between the two calculi are so strong and deep that we can address a logical problem by using transformations typical of algebra and interpret the result in logical terms.
Section 2 of the present paper is devoted to showing aspects of this correspondence\footnote{We are well aware that we here address a much debated and well known topic, still we think it is important to note certain similarities and/or dissimilarities}. 
Section 3 starts out from Boole's turning point: the recognition of a crucial difference between the two calculi, a difference that can be formalized in the \emph{law of duality}:
$$x x = x.$$

The product among classes is idempotent, which is not, as we know, among numbers. 
Boole assumes this law as the defining property of classes: if anything is a class, then the law of duality holds for it; on the algebraic side, $0$ and $1$, seen as numbers, fulfill the law of duality.
The correspondence can be re-established: the functions of this new (quantitative) algebra will be functions $f:\{0,1\}^n\to\IQ$. We call the resulting calculus, the \emph{pseudo-binary} calculus.  {This is in fact what Boole calculus really is: a pseudo-binary system, and not a purely binary one, as many may expect after the re-foundation of Boole's algebra pursued by Boole's successors!}

Now to a central question of Boole's project: given a function in the $\{\times, +, -, /, 0, 1\}$-language, what class does it represent ? The answer to this question is of fundamental importance and Boole answers it with \emph{the method of developments.} 
Let us start by considering the \emph{pseudo-binary} calculus limited to the 
functions $f:\{0,1\}^n\to \IZ$ in the language $\{\times, +, - ,0, 1\}$. We show that in this case the method of developments is equivalent to a variant of it that we call the \emph{method of intersections.} This variant we propose is very helpful to give a logical interpretation to the pseudo-binary calculus.

Just to get an intuitive idea, let us start by observing that given $n$ classes $x_1, \dots, x_n$, the universe can be partitioned into $2^n$ regions:

$ 1 = x_1 + (1 - x_1)$

$ 1 = x_1 x_2 + x_1(1-x_2) + (1-x_1) x_2 + (1-x_1)(1-x_2)$

$ 1 = x_1 x_2 x_3 + x_1x_2 (1-x_3) + x_1 (1-x_2) x_3 + (1-x_1) x_2 x_3 + x_1(1-x_2) (1-x_3)  +  (1-x_1) x_2 (1-x_3) + (1-x_1)(1-x_2) x_3 + (1-x_1)(1-x_2) (1-x_3)$

$\dots$

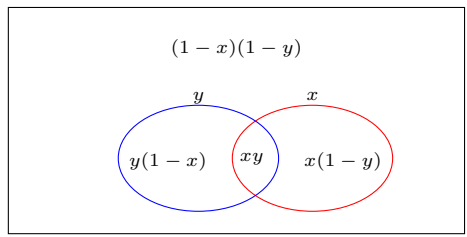
\begin{figure}[h]
\centering
\begin{tikzpicture}
\draw [color=red](05,01) ellipse (30pt and 20pt);    

\draw [color=blue] (03.5,01) ellipse (30pt and 20pt);    

\draw (01,00) rectangle (07,03);

\node at (5, 2) [below] {\scriptsize{$x$}};
\node at (3.5, 2) [below] {\scriptsize{$y$}};
\node at (04.2,01.2) [below] {\scriptsize{$xy$}};
\node at (03.1,01.2) [below] {\scriptsize{$y (1-x)$}};
\node at (05.4,01.2) [below] {\scriptsize{$x(1-y)$}};
\node at (04,02.7) [below] {\scriptsize{$(1-x)(1-y)$}};
	
\end{tikzpicture}
\caption{Development of 1 in two variables}
\label{fig:development_1}
\end{figure}

See Figure \ref{fig:development_1}. 

Then any arbitrary class $a = f(x,y)$ expressed by a function of the variables, say $x$ and $y$,  can be represented as a disjoint union as follows

$$ a x y +  a x(1-y) +  a (1-x) y + a (1-x)(1-y)$$

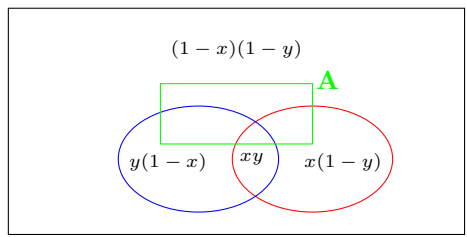
\begin{figure}[h]

\centering
\begin{tikzpicture}

\draw (01,00) rectangle (07,03);

\draw [color=red](05,01) ellipse (30pt and 20pt); 

\draw [color=blue] (03.5,01) ellipse (30pt and 20pt); 

\node at (04.2,01.2) [below] {\scriptsize{$xy$}};
\node at (03.1,01.2) [below] {\scriptsize{$y (1-x)$}};
\node at (05.4,01.2) [below] {\scriptsize{$x(1-y)$}};
\node at (04,02.7) [below] {\scriptsize{$(1-x)(1-y)$}};

\draw [color=green] (03,01.2) rectangle (05,02);    

\node [color=green] at (05.2,02.3) [below] {\bf A};


\end{tikzpicture}
\caption{The shadow of A}
\label{fig:shadow}
\end{figure}

This representation, given by the method of intersections,  is equivalent to what is obtained by the method of developments but at the same time it represents a class as the disjoint union of its "shadows" over the different regions into which the universe is partitioned.
See Figure \ref{fig:shadow}.

We immediately see that the set-theoretical union relevant here is the exclusive one. Boole has been severely blamed by Frege for having privileged the exclusive union instead of inclusive one, but Boole had very good reasons for his choice.

  {In Section 3.2 we show that if} we limit ourselves to consider functions  $f^n:\{0,1\}^n\to\{0,1\}$, we end up with nothing but the truth functions.  Now, if we interpret the variables on propositions, the product as $\wedge$, $+$ as $\vee$ and $1 - x$ as $\neg x$, then the development of $f^n(x_1, \dots, x_n)$ gives the 
proposition  in disjunctive normal form whose truth table is $f^n(x_1, \dots, x_n)$. In this case too, the use of exclusive disjunction is the natural one even if the inclusive one works as well.

There is more to the exclusive union and inclusive subtraction than a mere success in applying the method of intersections, these two operations are \emph{partial } operations, they are not defined when $x$ and $y$ are not disjoint or when $y$ is not a subset of $x$ and this, as we will see, corresponds perfectly to the cases in which the algebraic value of $x + y$ or $x - y$ is outside $\{0, 1\}$. The correspondence between the two calculi is amazing!  {In Section 3.3 we compare Boole's explanation of his strategy for treating non binary coefficients with our motivations based on the use of the methods of intersections.}

So far so good, but one could object that the logical calculus was intended to be a calculus of arbitrary classes, not just a calculus of the universal and the empty class only. It is well known from the standard propositional calculus and universal algebra that whatever can be said about arbitrary classes is reducible to what can be said about the universal and empty classes, see Section 3.2. Therefore as far as the logical principles are concerned, nothing is lost.

In  {Section 3.4} we comment on the very controversial solutions proposed by Boole when dealing 
with functions $f:\{0,1\}^n\to\IQ$ in the language also containing the operation of division. The main problem is how to interpret such terms as $0/0$ and $1/0$ possibly occurring in the developments of these functions. We show that these cases too admit of a logical interpretation by extending the equivalence between the method of developments and the method of intersections from an algebraic point of view.

\section{\emph{Logic and Algebra: a Formal Correspondence}}

In his celebrated work \emph{An Investigation of the Laws of Thought} [1854] George Boole analyzes the \emph{structural} similarities between quantitative algebra and logical reasoning in more {details} than he had {done} in his previous work \emph{The *Mathematical Analysis of Logic.}
His claim is the existence of evident \emph{formal} analogies shared by the universal laws of algebra and those of logic. Such formal analogies would be indeed so evident to be expressible in a common mathematical language in terms of equations having the same syntactical form. The semantic interpretation of any equation of this type will of course differ within the two different domains (quantitative algebra and logic); the revealed formal analogy does  not have in fact to correspond to some \emph{conceptual} identity. The \emph{syntactical} match is rather obtained in virtue of the use of the ordinary algebraic formalism as a unifying language between the two fields, and nothing can be said at the ontological level.  Boole himself openly expresses his empirical attitude of \emph{hypotheses non fingo} concerning the nature of this correspondence by pointing out that 

\begin{quote}...it is not affirmed that the process of multiplication in Algebra (...) possesses in itself any analogy with that process of logical combination wich $xy$ has been made to represent above; but only that if the arithmetical and logical process are expressed in the same manner, their symbolical expressions will be subject to the same formal law. (p.31) \end{quote} 

Same formal laws, different semantic interpretations: this is really the essence of the modern conception of formalism! The variables of the language will be interpreted either as numerical quantities, when our investigation concerns numerical algebra, or as classes of objects, when {our investigation is oriented} to logic, but the form of the respective laws will coincide.

\subsection {Logical operations and their interpretative problems}

The first example of correspondence between ordinary quantitative algebra and abstract logic that Boole shows concerns the multiplication symbol $\times$, whose logical interpretation turns out to correspond to the set theoretical operation of  \emph{class intersection} $\cap$\footnote{The legitimacy of an extensional approach to the treatment of concepts can be vindicated already from the beginning of Boole's \cite{B1854}:

\begin{quote}
(...) let $x$ represent ``all men'', or the class ``men''.  By a class is usually meant a collection of individuals, to each of which a particular name or description may be applied (p.28)
\end{quote}

His theory always admits in fact a set-theoretical interpretation and therefore can legitimately be viewed as a theory of classes. This does not mean that Boole really speaks the language of a contemporary set theorist. In this respect, the specific case of the ``class intersection'' operator even constitutes a particularly delicate case (see the following footnotes 5 and 6). It is nevertheless clear that we can immediately translate all the conceptual operations that he describes in the vocabulary of modern set theory, as we do in this case. We will adopt this approach throughout the paper, since we are mainly interested in investigating the value of Boole's system from a modern perspective.}:

\begin{quote}
Let it (...) be agreed, that by the combination $xy$ shall be represented that class of things to which the names or descriptions represented by $x$ and $y$ are simultaneously applicable.\footnote{As is customary in algebra, Boole always omits to write explicitly the occurrence of the symbol $\times$.}(p.28)
\end{quote} 

Boole claims that logical multiplication satisfies, \emph{formally}, the same \emph{law of commutativity} that holds of algebraic multiplication, or, more exactly, he affirmes that the two laws of commutativity holding in the two different domains can be expressed syntactically in the same  way, as $xy=yx$:

\begin{quote}
In the case of $x$ representing white things, and $y$ sheep, either of the members of this equation will represent the class of ``white sheep''. There may be a difference as to the order in which the conception is formed, but there is none as to the individual things which are comprehended under it. In like manner, if $x$ represent ``estuaries'', and $y$ ``rivers'', the expressions $xy$ and $yx$ will indifferently represent ``rivers that are estuaries'', or ``estuaries that are rivers'', the combination in this case being in ordinary language that of two substantives, instead of that of a substantive and an adjective as in the previous instance. Let there be a third symbol, as $z$, representing the class of things to which the term ``navigable'' is applicable, and any one of the following expressions,
\begin{center}
$zxy,zyx,xyz,$ ecc.
\end{center}
will represent the class of ``navigable rivers that are estuaries''. (p.29)
\end{quote}

Consistently with this, Boole observes that this law ``may be characterized by saying that the literal symbols $x,y,z$ are commutative, like the symbols of Algebra'' (p.31)\footnote{Boole is so used to the customary omission of the product symbol in the current algebraic notation that he apparently recognizes here no operation at all! He ascribes indeed the commutativity property not to the logical  operator in itself, but to the variables alone! Actually, although there might be really some ambiguity in this context, he refers explicitly to the \emph{logical process of combination} as compared to the \emph{arithmetical process of multiplication}  (p.31), as we have seen. Therefore, at least the existence of two analogue ``processes'' in the two respective fields is assumed; and if the arithmetical product is an operator, the same applies then to the logical combination.}

\subsubsection{Sum and subtraction as partial operations}

The quoted paragraph above is paradigmatic of Boole's peculiar {presentation} of his logical system: he does not introduce a ``static'' axiomatic calculus built on {basic} statements conceived as eternal \emph{a priori} truths. He rather investigates the rules of natural language and of the psychological construction of concepts almost as an \emph{experimental scientist}. The formal basic equations of logic are presented in the context of a colloquial and informal discussion which is mostly aimed at \emph{persuading} the reader about the plausibility of the results, and perhaps also at revealing the experience that led Boole to their discovery. For this reason, the main evidence for the fundamental logical laws that he formulates consists in concrete linguistic examples so obvious that nobody could raise any objection to their universal validity.

Hence, not surprisingly, the method applied to prove the commutativity of the product is directly extended to the sum:

\begin{quote}
We are not only capable of entertaining the conceptions of objects, [...] but also of forming the aggregate conceptions of a group of objects consisting of partial groups, each of which is separately named or described. For this purpose we use the conjunctions ``and'', ``or'', etc. ``Trees and minerals,'' ``barren mountains, and fertile vales,'' are examples of this kind. 
In strictness, the words ``and'', ``or'', interposed between the terms descriptive of two or more classes of objects, imply that those classes are quite distinct, so that no member of one is found in another. In this and in all other respects the words ``and'' ``or'' are analogous with the sign $+$ in algebra, and their laws are identical. Thus the expression ``men and women'' is [...] equivalent with the expression ``women and men''. Let $x$ represent ``men'', $y$, ``women''; and let $+$ stand for ``\emph{and'}' or ``\emph{or}'', then we have
$$x+y=y+x,$$
an equation which would equally hold true if $x$ and $y$ represented \emph{numbers}, and $+$ were the sign of arithmetical sum. (pp.32-33)
\end{quote}

One should not be confused by the use of ``and'' in this context: as is clear from the explanation, the operation {in question} is not class intersection, but rather the \emph{union} of (disjoint) classes. This is in fact consistent with an idiomatic use of the preposition ``and'' in our ordinary language, aimed at joining, although incorrectly, \emph{disjoint} notions.

We will call this particular set theoretical form of union ``\emph{disjoint union}''. As is well known, the preference accorded by Boole to the disjunctive interpretation of the union instead of its  \emph{inclusive} version (joining sets that may overlap) has met in the literature a lively criticism  {since Jevons' work (see} \cite{J74},  {p.70--71)}.  Even the founder of contemporary mathematical logic, G. Frege, defined Boole's choice a ``retrograde step away from Leibniz''.\footnote{Frege's \emph{Posthumous Writings}: ``On one point indeed Boole has taken a retrogade step away from Leibniz, in adding to the leibnizian meaning of $A + B$ the condition that the classes $A$ and $B$ should have no common element. W.Stanley Jevons, E. Schr\"oder and others have quite rightly not followed him in this.''( pag 10)}

In this paper we rather support the opposite view: arguments defending this and {others} of Boole's original {proposals} will be provided and discussed. It is true that in many circumstances his treatment of logic does not exactly coincide with the expectations of present day logicians, but this does not always entail a defect in its {pristine} design. We believe indeed that Boole should be regarded as ``something more'' than a \emph{clumsy} pioneer of the system {now called} ``Boolean algebra''. He  introduced a \emph{prolific} discipline whose goals are often different from those of our modern conception of logic, but with an independent value in themselves.

This holds for instance, in our opinion, of the disjoint union. We will not deny that this operation {apparently} carries a very unwelcome disadvantage: the occurrence of terms lacking of interpretations for some variable interpretations. An immediate example is ``$x+y$'', which is semantically undefined unless $x$ and $y$ denote disjoint classes. Nevertheless, there are several arguments in favour of Boole's choice.



Firstly, one could observe that symbolic logic is nowadays a subject that has gained a mature status of independence from traditional mathematics, whereas Boole's aim was that of showing that the connections linking logic and algebra were as close as they could be. In this particular respect, it is to be noted {in the first place} that the cardinality of a \emph{finite} set $A\cup B$ is equal to the sum of the cardinalities of $A$ and $B$ \emph{only if} $A$ and $B$ are disjoint sets.

Another fundamental similarity with ordinary algebra justifying the use of disjoint union can be found in its connection with its dual operator, introduced by Boole himself: the logical subtraction.

The latter constitutes in fact a sort of ``inclusive subtraction'', where the removed set must be completely included in the set from which it is substracted: the term $x-y$ will denote in fact the result of removing the class $y$ from the class $x$, under the assumption that $y$ is a subset of $x$. Again, the term ``$x-y$'' will be undefined unless this condition is fulfilled: 

\begin{quote}
(...) we cannot conceive it possible to collect parts into a whole, and not conceive it also possible to separate a part from a whole. This operation we express in common language by the sign \emph{except}, as, ``All men \emph{except} Asiatics,'' ``All states \emph{except} those which are monarchical.'' Here it is implied that the things excepted form a part of the things from which they are excepted. As we have expressed the operation of aggregation by the sign $+$, so we may express the negative operation above described by $-$ minus. Thus if $x$ be taken to represent men, and $y$, Asiatics, i. e. Asiatic men, then the conception of ``All men except Asiatics'' will be expressed by $x-y$. (pp.33-34)
\end{quote}

Still, these semantic restrictions make the logical subtraction be the formal inverse of the logical sum, like in integer algebra. Firstly, Boole observes that the transformation of $x=y+z$ into $x-z=y$ is valid:

\begin{quote}
Let us take the Proposition, ``The stars are the suns and the planets'', and let us represent stars by $x$, suns by $y$, and planets by $z$; we have then
$$x=y+z.$$
Now if it be true that the stars are the suns and the planets, it will follow that the stars, except the planets, are suns\footnote{Actually, as is clear from the context, Boole means not only that the remainder consists in \emph{some} suns, but, rather, in \emph{exactly all} suns.}. This would give the equation
$$x-z=y,$$
(...) Thus a term $z$ has been removed from one side of an equation to the other by changing its sign. This is in accordance with the algebraic rule of transposition. (pp.35--36)
\end{quote} 

The inference from $x=y+z$ to $x-z=y$ is valid in virtue of the disjointness of $y$ and $z$, as the example just mentioned shows (and Euler diagrams prove universally). 

In contrast, if $y$ and $z$ were allowed to overlap, the validity of this inference would no longer be deducible: by removing from $x$ all the elements which also are in $z$ it might be the case that only some, but \emph{not all}, members of $y$ are preserved. 

Hence, the restriction imposed to the logical sum justifies the given inference.

Secondly, let us consider the opposite direction, from $x - z = y$ to $x = y + z$. {Assume} $x - z = y$, if the class $z$ is entirely included in $x$, as Boole requires, by re-adding $z$ to the result of the previous subtraction, $y$, we indeed come back to $x$ (Euler-Venn diagrams can be used for a clarification). 

But  what would happen in modern set theory?

The definition of set subtraction used in this theory, as well as that of set union, avoids the occurrence of potentially vacuous terms: the set $A\setminus B$ is in fact defined as $A\cap \overline{B}$, where $overline{B}$ is the complementary set of $B$, hence it has always a denotation for every interpretation of $A$ and $B$. In particular, $B$ does not need to be contained in $A$.

Well, for this notion of subtraction, the {inverse} of Boole's inference would be no more justified: for $z$ possibly overflowing extensionally outside $x$, the result of adding $z$ to any given set, in particular to $y$ (i.e. $x-z$), may possibly exceed $x$.

The restriction imposed on the subtraction is then necessary to guarantee the {inverse} of Boole's inference.

In any case, the operator $\setminus$ coincides with Boole's subtraction under the same semantical restiriction, that is, when $B\sbsq A$.
\\

A further argument supporting Boole's conception of subtraction can be inspired by elementary arithmetic: within the realm of natural numbers, $x-y$ has in fact a denotation only when $y$ is a quantity entirely ``contained'' in the quantity $x$.
\\
\\
\subsubsection{Commutativity and associativity laws: a rift in the correspondence}


\begin{quote}
As it is indifferent for all the \emph{essential} purposes of reasoning whether we express excepted cases first or last in the order of speech, it is also indifferent in what order we write any series of terms, some of which are effected by the sign $-$. Thus we have, as in common algebra, $x-y=-y+x$. (p.34)
\end{quote}

In ordinary algebra the term $x-y$ can be transformed into $-y+x$ by applying the commutative law, but the validity of the equation $x-y=-y+x$ requires the introduction of negative numbers. We read in fact $-y+x$ as the sum of $-y$ with $+x$. We apply then the commutativity law for the \emph{sum} (and conventional syntactical abbreviations, such as that reducing, say, $(+a)+(-b)$ to {$(a-b)$} to {get} $x-y=(+x)+(-y)=(-y)+(+x)=-y+x$. Without the introduction of negative numbers and terms, the term $-y+x$ would simply be syntactically wrong, hence meaningless.

This fact is completely overlooked by Boole, who does not explain what a ``negative concept'' might be\footnote{The most intuitive solution, that is., interpreting negative concepts as ``complementary" classes, is not acceptable: Boole will define the complement of a class in a different manner, as we shall see later. For the moment, we only point out once again that the complementation operator is subject to no constraint, whereas a very strong constraint regulates the use of the subtraction sign, hence they must differ substantially.}, although, as we have seen, he sometimes explicitly refers to the \emph{signs} of the terms: ``a term $z$ has been removed from one side of an equation to the other by changing its sign'' (p.36), or ``any series of terms, some of which are effected by the sign $-$'' (p.34). 

To try and solve the dilemma, one may introduce, in virtue of a mere \emph{convention}, the term $-y+x$ as a syntactically well constructed expression denoting the same class as $x-y$. But this simple solution would hardly be applicable  to the general case of ``any series of terms, some of which are effected by the sign $-$'' (p.34). For example: what could be said about the term $x+y-z$? Would it have the same denotation as $x-z+y$? Actually, the term $x+y-z$ looks syntactically quite ambiguous: is the class $z$ removed from the class $x+y$ or from the class $y$ only? This ambiguity disappears in ordinary algebra, where the two terms $x+(y-z)$ and $(x+y)-z$ have the same denotation.  But what about Boole's logical calculus?

A little detour concerning the associative laws is in order, laws that Boole, to the best of our knowledge, never takes into consideration.
Actually, the applications of the associative laws in absence of occurrences of the subtraction symbol look definitely non problematic in the logical calculus, too. In fact $(x+y)+z=x+(y+z)$ is valid (for all suitable variable interpretations) as well as $(xy)z=x(yz)$ (always). The discrepancy between the two calculi really arises when negative terms {appear}, as with the case $x+y-z$ we started with. The equation $(x+y)-z=x+(y-z)$ is not universally valid in Boole's calculus, it fails in fact for those interpretations for which $z$ is a subset of $x+y$ but not of $y$ alone. Consequently, the search for a solution of our question whether

\begin{align}\label{star}
x+y-z = x-z+y
\end{align}

\noindent
should impose at least the use of explicit brackets.
\\
(a) Let us start with the reading of the first term as $x+(y-z )$. In this case $z$ is supposed to be a subset of $y$. Is there a reading of $x-z+y$ to justify (\ref{star})? The reading of $x-z+y$ as $(x-z)+y$ would in general be not allowed, since $z$ is not necessarily a subset of $x$. On the other hand, the reading as $x+(-z+y)$ even if potentially correct on the basis of the convention that $-z+y = y -z$ adds a new occurrence of $+$ which is not contained in the original term $x-z+y$, and moreover reads automatically $-z$ and $y$ as paired. 
This is admissible in algebra, but in the context of Boolean calculus it probably requires some \emph{ad hoc} syntactical convention.

(b) Consider now the reading of the first term as $(x+y)-z$. Is there a reading of $x-z+y$ to justify (\ref{star})? A variable interpretation according to which $z$ is a subset of $x+y$ but not of either $x$ or $y$ would be compatible with the given reading of the first term, but with none of the second, as both $(x-z)+y$ or $x+(-z+y)$ would not be acceptable.


One could then try and fix precise laws for an obligatory explicit use of brackets, but this solution would make Boole's calculus quite {clumsy} to handle. Boole ignores this and similar problematic aspects in the presentation of his system and instead concentrates, strategically, on those where a perfect correspondence between the two calculi can be exhibited, as in the case of the distributivity laws:

\begin{quote}
Let $x$ represent ``men'', $y$, ``women'' (...)

Let the symbol $z$ stand for the adjective ``European'', then since it is, in effect, the same thing to say ``European men and women'', as to say ``European men and European women'', we have
$$z(x+y)=zx+zy.$$
And this equation also would be equally true were $x$, $y$ and $z$ symbols of number (...) (p.33) 
\end{quote}

\begin{quote}
Still representing by $x$ the class ``men'', and by $y$ ``Asiatics'', let  $z$ represent the adjective ``white''. Now to apply the adjective ``white'' to the collection of men expressed by the phrase ``Men except Asiatics'' is the same as to say, ``White men, except white Asiatics''. Hence we have
$$z(x-y)=zx-zy.$$
This is also in accordance with the laws of ordinary algebra (p.34) 
\end{quote}

\section{The duality law and Boole's pseudo-binary calculus}

Boole's pivotal move is the recognition of a crucial difference between the two calculi, a difference that can be formalized in the \emph{law of duality}:
$$x x = x.$$

An explanation of this law is immediate: intersecting a set $x$ with itself results in the set $x$ itself. Hence, the product among classes is idempotent, which is not, as we know, among numbers. 

As {a} most striking discrepancy between the two domains,\footnote{ {See also} the equation (\ref{odot}) {as far as the division operator is concerned.}} this property will characterize very {efficiently} the notion of class in Boole's system: if anything is a class, then the law of duality will hold for it.\footnote{Boole apparently finds another break of the simmetry, when he observes the impossibility of infering, in his system, the truth of $x=y$ from that of $zx=zy$. Such an inference cannot indeed hold universally for all interpretations of $x$, $y$ and $z$ as classes. But after a little hesitation he is able to re-establish the desired analogy through the following remark: this inference is not even in algebra universally valid, although at a first superficial sight it might look to be so, in fact not acceptable for a vanishing $y$ (pp.36--37).}\label{impdiv}
%
%
%
%
%
%
%


But what follows now is a revolutionary idea with an enormous long lasting influence in the future development of logic and informatics\.footnote{Although Boole was  in this respect anticipated by Leibniz.}

It is well true that the law of duality does not hold, in general, over the whole domain of quantities, but, as he himself immediately points out, it can at least hold in some little fragment of it: more precisely on the sub-domain $\{0,1\}$! In this way, Boole is entitled to give a new life to his beloved correspondence: the match between the logical system and the ordinary algebraic calculus can be re-established \emph{provided that} the {interpretation} of the numerical variables {is} limited to $\{0,1\}$:

\begin{quote}
We have seen (...) that the symbols of Logic are subject to the special law,
$$x^2=x.$$
Now of the symbols of Number there are but two, viz. 0 and 1, which are subject to the same formal law. We know that $0^2=0$, and that $1^2=1$; and the equation $x^2=x$, considered as algebraic, has no other roots than 0 and 1. Hence, instead of determining the measure of formal agreement of the symbols of Logic with those of Number generally, it is more immediately suggested to us to compare them with the symbols of quantity \emph{admitting only of the value 0 and 1}. Let us conceive, then, of an Algebra in which the symbols $x,y,z$, etc. admit indifferently  of the values 0 and 1, and of these values alone. The laws, the axioms, and the processes, of such an Algebra will be identical in their whole extent with the laws, the axioms, and the processes of an Algebra of Logic. Difference of interpretation will alone divide them. Upon this principle the method of the following work is established. (pp.37-38)
\end{quote}

This passage is very important for two main reasons. 

(a) Boole introduces a new calculus, a sort of binary calculus that we will call ``\emph{pseudo-binary}''. It is true that variables can assume only the values 0 and 1, but this does not automatically extends to all terms. For instance, the term $x+y$ will assume the value 2 when both variables $x$ and $y$ are assigned the value 1. In this respect, Boole's logic (with operators $\times,+,-$)\footnote{In Section \ref{division} we will discuss the extension of the language to the division operator.} deals with functions $f:\{0,1\}^n\to\IZ$, rather then with functions $f:\{0,1\}^n\to\{0,1\}$.

{(b) The last statement sets the method that Boole will use throughout his book: the validity of any given logical law or logically correct  inference must be provable by a sequence of calculations in the new algebraic calculus, obtained from the usual one by adding instances of the duality law \emph{for all the variables}.}

{The origin of our modern conception of Boolean algebra as a logical binary calculus is to be found exactly in this strategy pursued by Boole!  
Still, this fact should not induce the reader to underestimate some irreducible differences between Boole's original view and our modern treatment of what is nowadays called Boolean algebra:}
(i) {Boole's pseudo-binary calculus has a quantitative interpretation, whereas  modern Boolean algebra is essentially qualitative and sees 0 and 1 not as numbers but rather as truth values (they can in fact be replaced by ``\emph{true}'' and ``\emph{false}''); logical operators are then free from formal connections with quantitative algebraic operations;\footnote{Although quantitative interpretations of Boolean operators are still possible, as $\hat\sigma(A\wedge B):= \min\{\sigma(A),\sigma(B)\}$, $\hat\sigma(A \vee B):=\max\{\sigma(A),\sigma(B)\}$ , $\hat\sigma(\neg A):=1-\sigma(A)$, for $\sigma(A),\sigma(B)\in\{0,1\}\sbsq\IN$.}}
(ii) it is well known that modern Boolean algebra has both a propositional and a set theoretical interpretation, Boole himself foresees the connections between the two fields, but he conceives this relation in a complete different way with respect to our current view, see Section 3.2.

We will be soon compelled to admit that the quantitative Boole's pseudo-binary calculus works really well, {after all.} A first motivation for this consists in the possibility of reducing it to a systems of equations of the ordinary algebraic calculus. More precisely, the assumption that all variables contained in a given equation $t=s$, say $x_1,...,x_n$, can be assigned only the values 0 and 1 can be expressed by the following system of equations of the ordinary algebra:
$$x_1^2=x_1$$
$$...$$
$$x_n^2=x_n$$
$$t=s.$$
Of course, the only possible solutions of this system assign to $x_1,...,x_n$ values in $\{0,1\}$. This is therefore a basic step to reduce Boole's pseudo-binary calculus to the ordinary algebraic one.
Although Boole never mentions it explicitly, this procedure is suggested by his own words ``the equation $x^2=x$, considered as algebraic, has no other roots than 0 and 1''. (p.37)\\

The double interpretation (quantitative and set-theoretical) proposed by Boole for his calculus does not concern only the binary operations, but also the constant symbols 0 and 1. He lets 0 denote the \emph{empty set}, and 1 the \emph{universal class}. Such interpretations are not only merely conventional, there are cogent motivations for that.  In particular, an immediate set-theoretical reading of the quantitatively valid equation $1\cdot x=x$ is that the intersection of the whole with any class $x$ results in the class $x$ itself. Beside the possible ontological justifications of interpreting 0 as the empty class (``no quantity=nothing=\emph{no thing}''), the equations $0\cdot x=0$ or $x+0=x$ will provide analogue meaningful set-theoretical motivations.

Since 1 is the \emph{whole}, $1-x$ will denote the complementary class of $x$.
The law of duality can now be re-written, via distributivity laws, as 
$$x(1-x)=0.$$ 
This equivalent notion expresses another fundamental property characterizing the notion of class: the intersection of a set with its complement results in the empty set. 

The law of duality must hold for every term denoting a class: it holds for variables \emph{in virtue of an assumption}, it holds for 0 and 1 \emph{numerically}. It holds as well for the complement of every class, as the next easy argument shows: let $t$ denote a class, then $t^2 = t$ by induction hypothesis, hence
 $(1-t)^2=(1-t)(1-t)=1-2t+t^2=1-2t+t=1-t$. Moreover, all terms denoting set intersections of classes, say $t_1 \cdot ... \cdot t_n$, satisfy the duality law: $(t_1 \cdot t_2 \cdot ... \cdot t_n)^2=t_1\cdot t_2 \cdot ... \cdot t_n\cdot t_1\cdot t_2 \cdot ... \cdot t_n=t_1^2 \cdot t_2^2 \cdot ... \cdot t_n^2=t_1\cdot t_2\cdot ... \cdot t_n$.
Consequenlty, the terms that Boole calls ``constituents'' satisfy the same law: given $n$ variables $v_1,...,v_n$, a \emph{constituent} is nothing but a product whose $i$-th factor is either $v_i$ or $(1-v_i)$.
For instance, all constituents expressible through $x$ and $y$ are $xy$, $x(1-y)$, $(1-x)y$, $(1-x)(1-y)$. 

More generally, every term $t$ representing a function $f:\{0,1\}^n\to\{0,1\}$ does denote a class, since $f(x_1,...,x_n)^2=f(x_1,...,x_n)$ holds necessarily. In other words, every expression that can assume only the values 0 and 1 for every possible interpretation $\sigma$ of the variables over $\{0,1\}$ will denote a class.\footnote{This does not mean that for any given term $t$ denoting a class each of its subterms will do the same. For instance $(1+1)-(1+1)$ satisfies the duality law (and can be viewed as a name for the empty set), nevertheless its subterm 1+1 does not. Of course, one may ask for an \emph{inductive} definition of the notion of \emph{class term}, in which case $(1+1)-(1+1)$ would no longer be admissible. But this requirement, which is not present in Boole's work, would restrict drastically the extension of the concept. A further example: $(x+x)-x$. This is essentially only a redundant re-writing of the basic term $x$, which represents a class for both binary interpretations of the  involved variable $x$. Nevertheless, the subterm $x+x$ satisfies the duality law only for the interpretation $\sigma(x)=0$. Our definition of class term as any term satisfying the duality law is more tolerant than the inductive one, and can be very naturally characterized in the following way. We call a term $t'$ \emph{primitive} if it can algebraically be no longer simplified. When a term $t$ satisfies the duality law (possibly with respect to a variable interpretation) and is reducible to the primitive term $t'$, then $t$ has the same denotation as $t'$ in virtue of the equation $t=t'$ (this convention is of course well defined in virtue of the laws of ordinary algebra, which imply the uniqueness of the primitive form $t'$ of $t$). Notice that all subterms in $t'$ will satisfy the duality law (possibly with respect to the considered variable interpretation). This approach looks very natural, but it will require some care for the case of the division operator.}

On the other hand, there exist terms whose values lie in $\{0,1\}$ only for some binary interpretations $\sigma$ of their variables. This is for example the case, as we know, of $x+y$ {but this is not a problem since the pseudo-binary calculus} is supposed {to go hand in hand with} the set theoretical {interpretation}: according to both, $x+y$ will indeed denote a class not always but only \emph{sometimes}. Even more remarkably, the correspondence between the two interpretations looks even stronger after a more careful inspection. On the one hand, $x+y$ does not satisfy the law of duality exactly for $\sigma(x)=\sigma(y)=1$, with 0 and 1 seen as numbers. Let us now have a deeper look at the set theoretical interpretation. The term $x+y$ denotes a set only when $x$ and $y$ are disjoint classes. The universal class is not disjoint with itself, whereas the empty set is disjoint with everything (even with itself!). Therefore, the sum of 1 with itself is even set theoretically not acceptable, while any sum involving the empty set as one of its addenda will be {automatically} admissible.\footnote{In the example just discussed, we find another justification for Boole's use of disjoint union: overlapping sets are not acceptable, just like $1+1=2\notin\{0,1\}$!}

This is a very appreciable property of Boole calculus: whatever can be said for arbitrary classes is reducible to what can be said for the empty and the universal sets alone! Some explanation for this apparently bizzarre ``set theoretical collapse'' is needed.\footnote{In contemporary algebra of logic this fact is a very well known consequence of the Stone's Representation Theorem. Still we want to give an intuitive motivation in terms of ordinary propositional calculus.}
 We know that in modern propositional calculus every formula admits of a double interpretation. According to the first option, variables are seen as \emph{arbitrary} propositions, and operation symbols $\wedge$, $\vee$, $\neg$ as propositional connectives. In the second option, the variables are seen as \emph{arbitrary} sets and the operation symbols are interpreted as the set operators ``union'' ($\cup$), ``intersection'' ($\cap$) and ``complementation'' ($\overline{\cdot}$). Nevertheless, one could mantain  that the first interpretation is somehow \emph{vacuously redundant}. In fact, via the usual binary truth value semantics, a propositional formula is reduced \emph{in reality} to a binary function over truth values. Hence, one could  \emph{get rid} of all that huge world of arbitrary propositions and \emph{keep only} a single tautology $\top$ and a single contradiction $\bot$. That is, the restriction of all possible variable interpretations to the domain $\{\top,\bot\}$ suffices to establish whether a given formula is a tautology, a contradiction, whether a formula is a logical consequence of others, and so on. For what concerns this paper, the following remark is of crucial importance: ordinary propositional logic and set theory define the same class of valid and invalid equations (where ``equation'' in the first case must of course be intended as ``semantic equivalence''\footnote{Every formula can obviously be reduced to a semantic equivalence via logical constants: for example, $A\vee\neg A$ can be read as $A\vee\neg A\equiv\top$, which corresponds set theoretically to $A\cup overline{A}=U$, for $U$ the universal set.}). Therefore, to prove the validity of the equation $A\cap (B\cup C)=(A\cap B)\cup(A\cap C)$, one can convert this formula into the propositional formula $A\wedge (B\vee C)=(A\wedge B)\vee(A\wedge C)$, and then check its validity only for those variable interpretations that assign to $A,B,C$  the statements $\top,\bot$. Correspondingly, the interpretation of $A,B,C$ as the universal set or the empty set must suffice to test the validity of $A\cap (B\cup C)=(A\cap B)\cup(A\cap C)$.

We can interpret Boole's operations as connectives whose meaning is fixed by the truth tables given below. For the connective $\times$ nothing new. Its binary interpretation coincides with the ordinary multiplication over $\{0,1\}$, and hence with the truth table of the propositional conjunction $\wedge$:
 
\vspace{0,2cm}

\begin{tabular}{p{0cm}p{0cm}p{0,5cm}p{1cm}}
$A$ & $B$ & & $A\times B$\\
1 & 1 & & 1\\
1 & 0 & & 0\\
0 & 1 & & 0\\
0 & 0 & & 0\\
\end{tabular}

\vspace{0,2cm}

This goes on a par with the fact that this operator corresponds to the ordinary set theoretical intersection $\cap$.

As to the operator $+$, its binary truth-table can only be \emph{partial} and this goes on a par with the fact that this operator corresponds to the disjoint union
\vspace{0,2cm}

\begin{tabular}{p{0cm}p{0cm}p{0,5cm}p{2cm}}
$A$ & $B$ & & $A+B$\\
1 & 1 & & {\small not allowed}\\
1 & 0 & & 1\\
0 & 1 & & 1\\
0 & 0 & & 0\\
\end{tabular}

\vspace{0,2cm}

First, observe that the arithmetically non allowed case  {reflects}, set theoretically, the unique case in which the  {corresponding} sets are not disjoint, and hence their union cannot be done.
Second, the table expresses a \emph{partial} truth function which coincides, over its partial domain, with the truth table of \emph{aut}, as well as  {of} \emph{vel}. And here is a curious fact. Although it has often been stated that the  {propositional interpretation of Boole's sum is the exclusive disjunction},\footnote{ {See for instance} \cite{S82} p.23, \cite{H91} p.12, \cite{B97},  {p.175.}} over its domain of definibility, it can be read \emph{indifferently} in both manners, as one prefers!\footnote{ {The acceptability of the inclusive interpretation has been pointed out already by J. Corcoran} (\cite{C86}  {p.71) and C. Badesa} (\cite{B04}  {p.3), who is even inclined to recognize it as the one really underlying Boole's operator. Actually, both versions are correct.}}

\vspace{0,2cm}

As to the operator $-$, its binary truth-table can again only be \emph{partial}, as well as its set theoretical interpretation as inclusive set subtraction:
\vspace{0,2cm}

\begin{tabular}{p{0cm}p{0cm}p{0,5cm}p{2cm}}
$A$ & $B$ & & $A-B$\\
1 & 1 & & $0$\\
1 & 0 & & 1\\
0 & 1 & & {\small not allowed}\\
0 & 0 & & 0\\
\end{tabular}

\vspace{0,2cm}
Observe that the arithmetically non allowed case corresponds, set theoretically, to the unique case in which the set to be removed is not contained in the other. 
\\

The truth table technique shows perfectly well that the empty and the universal set suffice for all semantic purposes, as they can simulate the behaviour of \emph{true} and \emph{false} with respect to the propositional interpretations of $\times,+,-$ and they are even able to distinguish between allowed and non allowed cases. 

This is enough for Boole's ``set theoretical collapse''. 

\subsection{The method of the developments}

Given an arbitrary functional term\footnote{Coherently with Boole's notation, we will not introduce a notational distinction between functions and functional symbols. The difference will be clear from the context.} defined in the $\{+,\times,-,1,0\}$-language, what class does it represent?
The method of the \emph{developments} transforms each functional term  in $n$ variables into a term in \emph{canonical form} i.e. a term expressed by the sum of $2^n$ products each of which is composed of a numerical value (a \emph{coefficient}) and  a \emph{constituent}.


Here is the development of $f(x_1, \dots, x_n)$:

\begin{align}\label{nabla}
f(x_1, \dots, x_n)=\sum_{1 \leq i \leq 2^n}f( \sigma^i(x_1) , \dots,  \sigma^i(x_n)) \cdot s_{i_1} \cdot \, \, ...  \cdot \, s_{i_n} 
\end{align}

%
for all $2^n$ possible binary functions $\sigma^i : \{ x_1, \dots, x_n\} \to \{0, 1\}$  and   $s_{i_k}, 1 \leq k \leq n,$ {where}:
\\
\\
$s_{i_k}= \left\{ \begin{array}{ll}
x_k & \quad
    \textrm{if} \quad \sigma^i(x_k)=1\\
    1 -x_k & \quad  \textrm{if} \quad
\sigma^i(x_k)=0
\end{array}  \right. $.
\\
\\
An example: consider $f(x_1, x_2)$; then
$$f(x_1, x_2)=\sum_{1 \leq i \leq 4}f(\sigma^i(x_1), \sigma^i(x_2)) \cdot s_{i_1} \cdot s_{i_2} =$$
$$ f( 1,1 ) \cdot x_1\cdot x_2 \, + \, f(1,0) \cdot x_1\cdot (1 - x_2) \, + \,   f(0,1) \cdot (1 - x_1)\cdot x_2  \, + \, f(0,0) \cdot (1 - x_1) \cdot (1 - x_2).
$$
\\
We have to show that the equality (\ref{nabla}) is valid, i.e.,  however we interpret the variables $x_1, \dots ,x_n$ that equality is true. This amounts to say that for every interpretation $\tau$ of the variables $x_1, \dots, x_n$ into $\{0,1\}$,

\begin{align}\label{nabla2}
f(\tau(x_1), \dots ,\tau(x_n))=\sum_{1 \leq i \leq 2^n}f( \sigma^i(x_1), \dots , \sigma^i(x_n)) \cdot \hat\tau(s_{i_1}) \cdot \, \, ...  \cdot \, \hat\tau(s_{i_n})
\end{align}

where $\hat\tau(x_k)=\tau(x_k)$, and $\hat\tau(1-x_k) = 1 - \tau(x_k)$. Notice that the $n$-tuple  $(\sigma^i(x_1), \dots , \sigma^i(x_n))$ is an $n$-tuple of 0 and 1, so the function $\tau$ plays no role on $f( \sigma^i(x_1), \dots , \sigma^i(x_n))$.

The function $\tau$ is necessarily equal to a $\sigma^j$ for exactly one $j$, $1 \leq j \leq 2^n$, so
\\
$f(\tau(x_1), \dots ,\tau(x_n)) = f(\sigma^j(x_1), \dots  ,\sigma^j(x_n)) $. Moreover we show that for that $j$: $\hat\tau(s_{j_1}) \cdot \, \,... \,  \cdot  \, \hat\tau(s_{j_n})$  is equal to 1. In fact for every $k$, $ 1 \leq k \leq n$, either $s_{j_k} = x_k $ or $s_{j_k} = 1 - x_k $.
In the first case $\sigma^j(x_k) = 1$, therefore $\hat\tau (s_{j_k} ) =1$;  in the second case, $ \sigma^j(x_k) = 0$, therefore $\hat\tau (s_{j_k} ) = 1 - \tau(x_k)=1 - \sigma^j(x_k) =1$.
\\
Now we show that for all $h \neq j$, $ 1 \leq h \leq 2^n$, $(\hat\tau(s_{h_1}) \cdot ... \cdot  \, \hat\tau(s_{h_n}))  = 0$.  For,  let $(\sigma^h(x_1), \dots , \sigma^h(x_n))\neq (\sigma^j(x_1) , \dots , \sigma^j(x_n))$, say $\sigma^h(x_k) \neq \sigma^j(x_k)$ for some $k$, and consider the constituent $s_{h_1}\cdot ... \cdot  \, s_{h_n}$. 
If $\sigma^h(x_k) =1$, then $s_{h_k} =x_k$ and $\sigma^j(x_k) = 0$ (since  $\sigma^h(x_k) \neq \sigma^j(x_k) $), therefore
$\hat\tau(s_{h_k})=\tau(x_k)$ = $\sigma^j(x_k) = 0$, and so the whole product $(\hat\tau(s_{h_1}) \cdot  ... \cdot  \, \hat\tau(s_{h_n}))$  is equal to 0.
Analogously, if $\sigma^h(x_k )=0$,  then $s_{h_k} = 1-x_k$ and $\sigma^j(x_k) = 1$ (since  $\sigma^h(x_k) \neq \sigma^j(x_k) $), therefore $\hat\tau(s_{h_k})=\hat\tau(1-x_k) = 1-\tau(x_k)=1-\sigma^j(x_k) = 0$, and so the whole product $\hat\tau(s_{h_1})\,  \cdot ... \cdot \, \hat\tau(s_{h_n})$ is equal to 0 again.
In conclusion, for every interpretation $\tau$ of the variables, the equation (\ref{nabla}) reduces to the  trivial identity 
$$f(\tau(x_1),\dots,\tau(x_n))=f(\sigma^j(x_1),\dots,\sigma^j(x_n)) \cdot 1$$ 
with $\tau(x_1)= \sigma^j(x_1),\dots,\tau(x_n)= \sigma^j(x_n)$. 
\\

Briefly, the given proof tells us that for any given interpretation $\tau$ the summation in  (\ref{nabla2})  reduces to that component whose coefficient is determined by the unique $\sigma^j$ equal to $\tau$.

If we use the notation $i_1, \dots,  i_n$ instead of $\sigma^i(x_1), \dots , \sigma^i(x_n)$, the equation just proved can be written as
\\
\begin{proposition}\label{developments}
$$f(x_1, \dots, x_n)=\sum_{i_1, \dots ,  i_n \in \{0,1\}}f( i_1, \dots , i_n) \cdot s_{i_1} \cdot ... \cdot \, s_{i_n}.$$
\end{proposition}
\noindent
The proof given by Boole in 
Chapter V of his book (pp.72-74) is not so general and uniform for all $n$, but we can say that our presentation captures its \emph{hidden} essence.
\\
\\
The point here is somehow delicate. Both Boole's original proof and our generalized version share a \emph{purely algebraic} nature and are indeed formulated within the pseudo-binary interpretation of the calculus. Therefore, one may also ask for an intuitive set theoretical characterization of the development technique for those functional terms $f(x_1,\dots,x_n)$ which allow a set theoretical interpretation\footnote{Also in this case, we will not be pedantic in distinguishing rigorously between terms  denoting sets and their corresponding denotations, which only are properly sets, or between functional symbols and the functions over the realm of sets they denote.}. A similar requirement is in this context even more relevant then in others because of the possible presence of  \emph{numerical coefficients} $f(i_1,\dots,i_n)\in\IZ$ possibly lying outside $\{0,1\}$, hence a set theoretical interpretation of the whole development becomes particularly problematic.

To try and accomplish this legitimate request, we will first of all prove, \emph{algebraically}, an interesting fact: the development of a term $f(x_1,\dots,x_n)$ \emph{coincides} with the {product} of that term with the development of 1 (in the same variables). By applying the distributive laws, we will even be able to show that the coincidence holds \emph{addendum by addendum}. The precise meaning of this claim needs a more detailed explanation, and, for the sake of simplicity, we limit our discussion to the case of two variables only (the generalization to any arbitrary number of variables can be then easily deduced). 
\\
By ordinary algebraic syntactical manipulations it can be shown that 

\begin{align}\label{box}
1=xy+x(1-y)+(1-x)y+(1-x)(1-y) 
\end{align}

This equation is valid for all possible values of $x$ and $y$, not necessarily in $\{0,1\}$, that is,  without any use of the duality law.
\\
By the ordinary distributivity laws we deduce then 
\\
$f(x,y) \,= \, f(x,y)\cdot 1=$
\\
$f(x,y)[xy+x(1-y)+(1-x)y+(1-x)(1-y)]=$
\\
$f(x,y)xy+f(x,y)x(1-y)+f(x,y)(1-x)y+f(x,y)(1-x)(1-y)$. 
\\
\\
On the other hand, by the method of the developments 
\\
$f(x,y)=f(1,1)xy+f(1,0)x(1-y)+f(0,1)(1-x)y+f(0,0)(1-x)(1-y)$. 
\\
\\
Therefore, the two sums

 \begin{align}\label{clubsuit}
f(x,y)xy+f(x,y)x(1-y)+f(x,y)(1-x)y+f(x,y)(1-x)(1-y)
\end{align}

and 

 \begin{align}\label{clubsuit_clubsuit}
f(1,1)xy+f(1,0)x(1-y)+f(0,1)(1-x)y+f(0,0)(1-x)(1-y)
\end{align}

must denote, \emph{numerically}, the same quantity for each \emph{binary} interpretation $\sigma$ of the variables. We will prove something more, i.e., that \emph{each addendum of the first sum can be transformed, by algebraic calculations, into the corresponding addendum of the second sum} (and vice versa).

Before proving it, we want to give a set theoretical interpretation of this crucial result. To this purpose, we notice first of all that the right term of (\ref{box}) coincides with the \emph{development} of 1, since for arbitrary $i$ and $j$ we have $f(i,j)=1$ where $f$ is the constant function 1.\footnote{For each $n\geq 1$, the coefficient 1 can be identified with the $n$-ary function $f(x_1, \dots , x_n)$ whose value is constantly 1.} 

From the set theoretical point of view, the development of the universal class 1 can be interpreted as its decomposition into four disjoint classes, whose disjoint union results  indeed in the universal class itself. This fact can be easily checked by ordinary Eulero-Venn diagrams. See figure \ref{fig:development_1}. 
%
%
%
%
%
As we know, the set theoretical counterpart of the algebraic product is the set intersection. Therefore, each addendum of (\ref{clubsuit}) can be seen as the \emph{intersection of the class $f(x,y)$ with the corresponding region of the universal set} (whenever $f(x,y)$ is a class). The final sum is of course $f(x,y)$, algebraically, but this is perfectly coherent with the set theoretical interpretation. Let $f(x,y)$ express the class $A$,  if we relativize all regions of the universal set to $A$ (i.e. we keep only those portions overlapping with $A$), and then we join them, we obtain exactly the set $A$ itself. Again, the Eulero-Venn diagrams give a visual image of this fact. See figure \ref{fig:shadow}.
%
%
%
%
%
%
%
%
%
%
%
%
%

In virtue of this set theoretical interpretation, we call the method of analyzing $f(x,y)$ as specified by (\ref{clubsuit})  the \emph{method of the intersections} and we prove that

\begin{proposition}\label{developments-intersections-equivalence-a}
\emph{For functional terms $f(x,y)$ in the $\{+,\times,-,0,1\}$-language, the} method of the developments \emph{and the} method of the intersections \emph{are equivalent.}
\end{proposition}
This methodological equivalence makes the method of the developments more intuitive in the domain of sets: \emph{calculating the development of a set $f(x,y)$ coincides with re-constructing that set by joining disjoint ``pieces'' of it, each piece being the ``shadow'' of $f(x,y)$ over one of the four parts into which the universal set is standardly decomposed}\footnote{Of course, properly speaking, such four parts change along with the interpretation of $x$ and $y$.}.

Now to the proof of our claim. We simplify the notation for constituents:  we will write $c_{11}, c_{10}, c_{01}, c_{00} $ for $xy$, $x(1-y)$, $(1-x)y$ and $(1-x)(1-y)$, respectively. Whence, as announced, we prove that $f(xy)c_{ij}$ and $f(i,j)c_{ij}$ coincide up to algebraic trasformations\footnote{This means that they share the same numerical denotations, therefore we show each step of their reciprocal trasformation by using the identity symbol rather than some artificial symbol of syntactical reduction:  to manipulate the terms, we use nothing but algebraic \emph{equations}.} for $i,j \in \{0,1\}$ ($f(i,j)$ not necessarily in $\{0,1\}$).

The proof is then by induction on the syntactical construction of $f$:

\begin{itemize}
\item  $f(x,y):=0$: then $f(xy)c_{ij}=0=f(i,j)c_{ij}$; 
\item  $f(x,y):=k\neq 0$: then $f(xy)c_{ij}=kc_{ij}=f(i,j)c_{ij}$; 
\item $f(x,y)$ is a projection, for instance $f(x,y):=x$. If $c_{ij}$ is of the type $x\hat y$ for $\hat y\in\{y,(1-y)\}$, then it is associated with the coefficient $f(1,j)=1$. Therefore $f(i,j)c_{ij}=c_{ij}$; on the other hand it holds $f(x,y)c_{ij}=xc_{ij}=x^2\hat y= x\hat y=c_{ij}$. If $c_{ij}$ is of the type $(1-x)\hat y$ for $\hat y\in\{y,(1-y)\}$, then it is associated with the coefficient $f(0,j)=0$. Therefore $f(i,j)c_{ij}=0$; on the other hand it holds $f(x,y)c_{ij}=xc_{ij}=x(1-x)\hat y=0\cdot \hat y=0$;
\item $f(x,y):=f_1(x,y)\cdot f_2(x,y)$. Of course $f(i,j)=f_1(i,j)\cdot f_2(i,j)$. By IH $f_l(x,y)c_{ij}= f_l(i,j)c_{ij}$ for $l=1,2$. 
Then, by commutativity of the product and by the fact that $c_{ij}$ satisfies the duality law, $f(x,y)c_{ij}=(f_1(x,y)\cdot f_2(x,y))c_{ij}=(f_1(x,y)\cdot f_2(x,y))c_{ij}^2= f_1(x,y)c_{ij}\cdot f_2(x,y)c_{ij}=
f_1(i,j)c_{ij}\cdot f_2(i,j)c_{ij}=(f_1(i,j)\cdot f_2(i,j))c_{ij}^2=(f_1(i,j)\cdot f_2(i,j))c_{ij}= f(i,j)c_{ij}$;
\item $f(x,y):=f_1(x,y)\diamond f_2(x,y)$ for $\diamond\in\{+,-\}$. Of course $f(i,j)=f_1(i,j)\diamond f_2(i,j)$. By IH $f_l(x,y)c_{ij}= f_l(i,j)c_{ij}$ for $l=1,2$. Then, by distributivity laws, $f(x,y)c_{ij}=(f_1(x,y)\diamond f_2(x,y))c_{ij}= f_1(x,y)c_{ij}\diamond f_2(x,y)c_{ij}= f_1(i,j)c_{ij}\diamond f_2(i,j)c_{ij}= (f_1(i,j)\diamond f_2(i,j))c_{ij}=f(i,j)c_{ij}$.
\end{itemize}

The proof of the coincidence between the two methods is proved for the general case of generic functions that are not necessarily logical. When non logical functions are involved (for instance $f(x,y)=5$, $f(x,y):=-x-y$ or $f(x,y):=x+y+2$)  products will be seen as purely algebraic operations and not as set intersections. We nevertheless use in every situation the expression ``method of the intersections'', as the proved algebraic result represents indeed a generalization of the set theoretical perspective.

\vspace{0,2cm}

Let us now instantiate some specific case. In the following calculations, which show the equalence of the terms $f(x,y)c_{ij}$ and $f(i,j)c_{ij}$ for some concrete examples, we may even go beyond the rigid sequence of steps prescribed by the induction scheme of the proof of Proposition \ref{developments-intersections-equivalence-a}. After all, that induction scheme shows that $f(x,y)c_{i_j}$ and $f(i,j)c_{ij}$ denote, numerically, the same quantity; hence we can use all admissible algebraic transformations \emph{freely}:

\begin{itemize}
\item $f(x,y):=xy$ is always a class. It holds 
$$f(x,y)c_{11}=xyc_{11}=xyxy=xy=1\cdot xy=1\cdot c_{11}=f(1,1)c_{11}$$
$$f(x,y)c_{10}=xyc_{10}=xyx(1-y)=0=0\cdot c_{10}=f(1,0)c_{10}$$
$$f(x,y)c_{01}=xyc_{01}=xy(1-x)y=0=0\cdot c_{01}=f(0,1)c_{01}$$
$$f(x,y)c_{00}=xyc_{00}=xy(1-x)(1-y)=0=0\cdot c_{00}=f(0,0)c_{00}$$

\item $f(x,y):=x+y$ is sometimes a class. It holds anyway
$$f(x,y)c_{11}=(x+y)xy=x^2y+xy^2=xy+xy=2xy=f(1,1)c_{11}$$
$$f(x,y)c_{10}=(x+y)x(1-y)=(x+y)(x-xy)=x^2-x^2y+xy-xy^2=$$ $$=x-xy+xy-xy=x-xy
=x(1-y)=1\cdot x(1-y)=f(1,0)c_{10}$$
\\
$$f(x,y)c_{01}=(x+y)(1-x)y=(x+y)(y-xy)=xy-x^2y+y^2-xy^2= $$ $$= xy-xy+y-xy=y-xy=$$
$$=(1-x)y=1\cdot (1-x)y=f(0,1)c_{01}$$
\\
$$f(x,y)c_{00}=(x+y)(1-x)(1-y)=(x+y)(1-y-x+xy)= $$ $$=x-xy-x^2+x^2y+y-y^2-xy+xy^2= $$ $$=x-xy-x+xy+y-y-xy+xy=0=0\cdot(1-x)(1-y)=f(0,0)c_{00}$$
\\
\end{itemize}


\subsection{Boole's propositional calculus}\label{Booleprop}

Let us now draw our attention to the cases in which $f(x_1,\dots,x_n)$ \emph{always} satisfies the law of duality, that is, for every possible interpretation of the variables; this implies that all coefficients of type $f(i_1, \dots, i_n)$, for $i_1, \dots , i_n\in\{0,1\}$, belong to $\{0,1\}$. In this context {a} set theoretical interpretation of the whole development of $f(x_1,\dots,x_n)$ is {immediate}. Let $f(i_1,...,i_n) \cdot c_{i_1...i_n}$ be the $i$-th element of the development. By assumption, $f(i_1,...,i_n)$ must be interpreted {either} as the universe or the empty set. Correspondingly, its intersersection with the consituent  $c_{i_1...i_n}$ will result {either in the constituent itself or in the empty set}. In other words, the \emph{coefficient} $f(i_1,...,i_n)$ will determine whether the corresponding constituent $c_{i_1...i_n}$ will contribute or not to the construction of the set $f(x_1,...,x_n)$ as a member of the development union: $f(x_1,...,x_n)$ will be the set resulting from the union of the {``survived'' }constituents.



In such cases, the calculation of developments reflects exactly the standard method of constructing a propositional formula corresponding to a given truth table. As well known, such formula is in disjunctive normal form. Take the $i$-th line of the truth table with value 1, construct the conjunction of $n$ literals $l_1,\dots,l_n$, with $l_j\equiv A_j$ if in that line $A_j$ is assigned the value 1, $l_j\equiv \neg A_j$ if it is assigned the value 0. Build the disjunction of all the conjuncts so obtained.

\vspace{0,2cm}

For example, a formula corresponding to the truth table

\vspace{0,2cm}

\begin{tabular}{p{0,2cm}p{0,2cm}p{0,2cm}p{0,5cm}p{0,2cm}}
$A_1$ & $A_2$ & $A_3$ & \\
1 & 1& 1& & 1\\
1 & 1& 0& & 0\\
1 & 0& 1& & 1\\
1 & 0& 0& & 1\\
0 & 1& 1& & 0\\
0 & 1& 0& & 0\\
0 & 0& 1& & 1\\
0 & 0& 0& & 1\\
\end{tabular}

\vspace{0,2cm}

is the following: $(A_1\wedge A_2\wedge A_3)\vee(A_1\wedge\neg A_2\wedge A_3)\vee(A_1\wedge\neg A_2\wedge\neg A_3)\vee(\neg A_1\wedge\neg A_2\wedge A_3)\vee(\neg A_1\wedge\neg A_2\wedge\neg A_3)$. In Boole's language, every positive literal $A_j$ is replaced by $x_j$ and every negative literal $\neg A_j$ is replaced by $(1-x_j)$, $\wedge$ is replaced by $\times$ and $\vee$ by $+$. Therefore, if $f(x,y,z)$ is the (binary) function whose graph coincides with the given truth table, we obtain that $f(x,y,z)=xyz+x(1-y)z+x(1-y)(1-z)+(1-x)(1-y)z+(1-x)(1-y)(1-z)$. This is equal to $1\cdot xyz+0\cdot xy(1-z)+1\cdot x(1-y)z+1\cdot x(1-y)(1-z)+0\cdot(1-x)yz+0\cdot(1-x)y(1-z)+1\cdot(1-x)(1-y)z+1\cdot(1-x)(1-y)(1-z)$, that is, the development of $f(x,y,z)$.\footnote{One could argue that in this way $\vee$ is translated into the exclusive disjunction rather than into the inclusive one, but this is irrelevant: the different lines of the truth table are mutually exclusive, hence in this context ``vel'' and ``out'' coincide!}
By interpreting set theoretically the development so obtained, according to the instructions suggested above, we have {automatically} a set theoretical interpretation of $f(x,y)$ as a well determined set.

Vice versa, every development with all coefficients $f(i_1,...,i_n)$ in $\{0,1\}$ defines a truth table by retracing the process in the opposite direction.

{Boole comes very close to} a fundamental aspect of the {symmetry between} propositions and classes, at the base of the modern Boolean algebra.


Nevertheless, he does not seem to be aware of this fact. The correspondence between the two domains that he {sees} is of completely different nature. Well, we could even argue that he sees no correspondence at all! Rather, only a trivial reduction; the idea is to interpret a variable $v$ as the set of moments in which a given proposition $A$ is true:

\begin{quote}
Let us employ the capital letters $X, Y, Z$, to denote the elementary propositions concerning which we desire to make some assertion touching their truth or falsehood (...) And let us employ the corresponding small letters $x,y,z$, considered as expressive of mental operations, in the following sense, viz. : Let $x$ represent an act of the mind by which we fix our regard upon that portion of time for which the proposition $X$ is true; and let this meaning be understood when it is asserted $x$ \emph{denotes} the time for which the proposition $X$ is true. Let us further employ the connecting signs $+,-,=, \&$ in the following sense, viz.: Let $x+y$ denote the aggregate of those portions of time  for which the propositions $X$ and $Y$ are respectively true, those times being entirely separated from each other. Similarly let $x-y$ denote the remainder of time which is left when we take away from the portion of time for which $X$ is true, that (by supposition) included portion for which $Y$ is true. Also, let $x=y$ denote that the time for which the proposition $X$ is true, is identical with the time for which the proposition $Y$ is true. (pp.164--165)
\end{quote}

\subsection{Empty constituents}

What does it happen when the coefficients $f(i,j)$ {do} not belong to $\{0,1\}$? 

For instance, the \emph{partial} binary tables of $+$ and $-$ presented earlier can be extended to \emph{total} tables by completing the missing arithmetical calculations, so to obtain

\vspace{0,2cm}

\begin{tabular}{p{0cm}p{0cm}p{0,5cm}p{1cm}p{1cm}p{1cm}}
$A$ & $B$ & & $A+B$ & & $A-B$\\
1 & 1 & & $2$ & & $0$\\
1 & 0 & & 1 & & 1\\
0 & 1 & & 1& & -1\\
0 & 0 & & 0 & & 0\\
\end{tabular}

\vspace{0,2cm}

This point is very delicate, since the values 2 and $-1$ have no direct set theoretical meaning (nothing is greater than the whole or smaller than the nothing!). {Contrary to expectations, Boole is able to draw from such ``unlogical'' values, a kind of} \emph{logical information}. The idea is {simple}: every constituent {preceded} by a coefficient different from 0 or 1 in the development of a term {is going to be} identified with the empty set! Be careful now: it is the constituent, and not its coefficent, to be identified with such a set. The, so to say, set theoretical interpretation of the coefficient reduces to a mere \emph{warning} for the {emptiness} of its constituent. 

Let us calculate for example the development of the term $f(x,y):=x+y$:
$$x+y=2xy+x(1-y)+(1-x)y.$$
The constituent $xy$ has coefficient 2. According to Boole's claim, this coefficient indicates that the intersection of $x$ with $y$ is the empty set (as really must be, being $x$ and $y$ by assumption disjoint!). As another example, consider the development of $x-y$. This results in
$$x-y=x(1-y)-(1-x)y.$$
Here, the intersection of $(1-x)$ with $y$ will be again the {empty set} (as must be, by the assumption that the whole $y$ is included in $x$).

In spite of its good behaviour in these two {examples}, the reader might feel quite perplex about a suggestion that looks totally arbitrary in itself. But this is not the case, and Boole justifies his claim by proving, algebraically, the following theorem {which opens} the 11th paragraph of the sixth chapter:

\vspace{0,5cm}

\emph{If a function $V$, intended to represent any class or collection of objects, $w$, be expanded, and if the numerical coefficient, $a$, of any constituent in its development, do not satisfy the law, $a(1-a)=0$, then the constituent in question must be made equal to 0.}

\vspace{0,5cm}

In the proof, Boole assumes that

\begin{align}\label{star}
w=a_1t_1+a_2t_2+ \dots +a_nt_n
\end{align}


where the right term of this equation is the development of $V$ (which means that the terms $t_i$ are constituents).  From (\ref{star}) we immediately obtain that $w^2=(a_1t_1+a_2t_2+...+a_nt_n)^2$,
and then $w=(a_1t_1+a_2t_2+...+a_nt_n)^2$, because $w$, as a variable, is subjected to the duality law.
The distributivity laws allow us to compute the value of the square, i.e., of $(a_1t_1+a_2t_2+...+a_nt_n)(a_1t_1+a_2t_2+...+a_nt_n)$. The computation is very easy because, on the one hand, $a_it_i \cdot a_it_i = a_i^2 t_i$ for all $i$, and on the other,  $a_it_i \cdot a_jt_j = 0$ whenever $i \neq j$: for that, observe that for some variable $v$, either $v$ is a factor of $t_i$ and $(1-v)$ of $t_j$, or viceversa.
\\
So we obtain
\begin{align}\label{star_star}
w=a_1^2t_1+a_2^2t_2+...+a_n^2t_n
\end{align}

Let $a_{i_1},...,a_{i_m}$ be all the coefficients among $a_{1},...,a_{n}$ that do not satisfy the duality law. For all other coefficients it holds $a=a^2$, therefore by subtracting (\ref{star_star}) from (\ref{star}) we deduce 
\\
\\
$0=(a_{i_1}-a_{i_1}^2)t_{i_1}+(a_{i_2}-a_{i_2}^2)t_{i_2}+...+(a_{i_m}-a_{i_m}^2)t_{i_m}$, 
that is to say 

\begin{align}\label{star_star_star}
0=a_{i_1}(1-a_{i_1})t_{i_1}+a_{i_2}(1-a_{i_2})t_{i_2}+...+a_{i_m}(1-a_{i_m})t_{i_m} 
\end{align}

By multiplying both sides in (\ref{star_star_star}) by $t_{i_j}$ with $1\leq j\leq m$, one infers in each case $0=a_{i_j}(1-a_{i_j})t_{i_j}$ (all the other addenda disappear for the very same reason as above). But $a_{i_j}(1-a_{i_j})\neq 0$ by assumption, and this yelds that $t_{i_j}=0$.


Boole's proof is algebraically correct. Notwithstanding, the reader may feel unhappy with it, as, once again, its purely algebraic nature hides completely any set theoretical interpretation. 

We give two alternative explanations closely related to each other which are  still of numerical nature,\footnote{It could not be any different, as long as coefficients with no direct set theoretical interpretation are possibly concerned.} nevertheless they are based on the coincidence of the method of the developments with that of the intersections. We think that this may help, in particular in the second argument, to \emph{perceive} some  {more evident}
set theoretical \emph{flavour}.  Again, for the sake of simplicity, we deal with the case of two variables only.

By (the proof of) Proposition \ref{developments-intersections-equivalence-a}, $f(i,j)c_{ij}=f(x,y)c_{ij}$ (for all possible binary variable interpretations). Morevorer we know that Boole's calculus can be simulated within the ordinary algebraic calculus by considering only the variable interpretations  which are solutions of the equation system $\bigstar$: 


\begin{eqnarray*}
x^2 &=& x\\
y^2 &=& y\\
f(x,y)^2 &=& f(x,y).
\end{eqnarray*}

First argument. For all binary variable interpretations satisfying $\bigstar$
it holds $f(x,y)^2=f(xy)$. For such interpretations, we deduce then:

\noindent
$f(i,j)c_{ij}=f(x,y)c_{ij}=f(x,y)^2c_{ij}^2=(f(x,y)c_{ij})^2=(f(i,j)c_{ij})^2$. 

Therefore $f(i,j)c_{ij}=(f(i,j)c_{ij})^2$. If $f(i,j)\neq f(i,j)^2$, then it must be $c_{ij}=0$ (for all binary interpretations allowed by  $\bigstar$.


%
%


Second argument. For all variable interpretations fulfilling  $\bigstar$, we can see $f(x,y)$ as a set. Hence, $f(x,y)c_{ij}$ will represent the intersection of two sets, which in turn is a set. Therefore, $f(x,y)c_{ij}$ must satisfy the duality law, and the same must then apply, numerically, to $f(i,j)c_{ij}$. When $f(i,j)\notin\{0,1\}$, the only condition for this to hold is that $c_{ij}$ is algebraically null, i.e., set theoretically seen as the empty set.


\subsection{The logical meaning of division}\label{division}

%
%

According to Michael Dummett, the division operator has \emph{no logical meaning} at all, in \cite{D59} Dummett writes ``He introduced a division sign for the operation inverse to intersection, and never succeeded in unravelling the complicated tangles which resulted from this.''

 Let us consider Boole's analysis of the proposition ``Clean beasts are those which both divide the hoof and chew the cud'', formalized as $x=yz$ (pp.86--87). Boole is tempted to infer from this equation the validity of a new one, $z=\frac xy$, through the ordinary algebraic introduction of division; but immediately after he observes that this ``equation is not at present in an interpretable form'' (p.87). At the same time, {he goes on saying} ``If we can reduce it to such a form it will furnish the relation required'' (p.87). 
%
%
\\
The interpretable form Boole is referring to is the development of $\frac xy$, which is:

\begin{align}\label{paragrafo}
\frac xy \, = \, \frac11xy+\frac10x(1-y)+\frac01(1-x)y+\frac00(1-x)(1-y)
\end{align}


By that, the problem of providing a meaning to the operation of division is reduced to that of providing a meaning to $\frac 1 0$ and $\frac 00$.

\

According to Boole the quotient $\frac10$ should be dealt with in the same way as non binary coefficients seen so far. While commenting on his theorem for the treatment of non binary coefficients, he states:


\begin{quote}
(...) it may be shown generally that any constituent whose coefficient is not subject to the same fundamental law as the symbols themselves\footnote{Of course, the duality law.} must be separately equated to 0. The usual form under which such coefficients occur is $\frac10$. This is the algebraic symbol of infinity. Now the nearer any number approaches to infinity (allowing such an expression), the more does it depart from the condition of satisfying the fundamental law above referred to. (p.91)
\end{quote}

Boole treats the ``\emph{infinity}'', as he calls it, as the limit of quantities whose degree of compliance of the duality law is proportionally inverse to their magnitude  (``the nearer any number approaches to infinity... the more does it depart from the condition of satisfying the fundamental law'':  formally, $\lim_{x\to\infty}(x^2-x)=\infty$). Hence, whatever this object might be, it should not satisfy the duality law at the highest degree; by that, its treatment as an \emph{ordinary} number different from 0 and 1 will follow.\footnote{In ordinary algebra, the argument $(\frac10)^2=\frac{1^2}{0^2}=\frac10$, \emph{apparently} proving the duality law for $\frac10$, has no sense, as $\frac10$ is no number. Consequently, the argument does not need to hold in Boole's calculus either. One can also observe that in non standard analysis the sequences $\Big(\frac1{\frac1n}\Big)_n$ and $\Big(\frac1{(\frac1n)^2}\Big)_n$ define two \emph{well distinct} unlimited non standard numbers, the second being larger than the first one.}


As for the coefficient $\frac00$ Boole affirms:

\begin{quote}
The symbol $\frac00$ [...] does not necessarily disobey the law we are here considering, for it admits of the numerical values 0 and 1 indifferently. Its actual interpretation, however, as an indefinite class symbol, cannot, I conceive, except upon the ground of analogy, be deduced from its arithmetical properties, but must be established experimentally. (pp.91--92)
\end{quote}

\begin{quote}
[...]  The symbol $\frac00$ indicates that a perfectly \emph{indefinite} portion of the class, {i.e.} \emph{some}, \emph{none}, or \emph{all} of its members are to be taken. (p.92)
\end{quote}


Boole claims that the interpretation of $\frac 0 0$ ``must be established experimentally'' (p.92), and starts by considering the statement ``Men who are not mortals do not exist''. By denoting ``men'' as $y$ and ``mortal beings'' as $x$ he comes to the equation $y(1-x)=0$, transformed into $y-yx=0$,  and finally into $yx=y$. In ordinary algebra one would then proceed by dividing both sides by $y$ (for a non vanishing $y$), so to obtain $x = \frac y y $ and conclude $x = 1$. But, of course, the conclusion that all things are mortal is not contained in the assumption and we do not want it!

Boole knows that ``the operation of division cannot be \emph{performed} with the symbols with which we are now engaged'' (p.89)
and again he suggests us to calculate the development of $\frac yy$: ``Our resource, then, is to \emph{express} the operation, and develop the result by the method' of the preceding chapter'' (i.e. the method of developments) (p.89). 

%

Well, from the equation $x=\frac yy$ by calculating the development of $\frac yy$ we obtain:
$$x=y+\frac00(1-y).$$
Simple semantic observations lead Boole to extract some meaningful logical information from it:

\begin{quote}
This implies that mortals $(x)$ consist of all men $(y)$, together with such a remainder of beings which are not men $(1-y)$, as will be indicated by the coefficient $\frac00$. Now let us inquire what remainder of ``not men'' is implied by the premise. It might happen that the remainder included all the beings who are not men, or it might include only some of them, and not others, or it might include none, and any one of these assumptions would be in perfect accordance with our premiss [...] and therefore the expression $\frac00$ here indicates that \emph{all}, \emph{some} or \emph{none} of the class to whose expression it is affixed it must be taken. (pp.89--90)
\end{quote}

He then quite optimistically and without proof concludes:

\begin{quote}
Although the above determination of the significance of the symbol $\frac00$ is founded only upon the examination of a particular case, yet the principle involved in the demonstration is general, and there are no circumstances under which the symbol can present itself to which the same mode of analysis is inapplicable.
We may properly term $\frac00$ an \emph{indefinite class symbol}, and may, if convenience should require, replace it by an uncompounded symbol $v$, subject to the fundamental law, $v(1-v)=0$. (p.90)
\end{quote}

The task of providing a clear meaning to the operation of division so as to include also expressions such as $\frac 0 0 $ and  $\frac 1 0 $ remains untouched by Boole.  In what follows we try to address this task and we make a proposal with the effect that the `new' operation of division will coincide with the standard operation whenever the denominator is different from 0, and it is in accordance with Boole's \emph{desiderata} as to $\frac k 0 $, $k \geq 0$.

 Our {proposal} starts from the obvious consideration that the division operator is the inverse operation of the multiplication. This is essentially what Boole himself requires, in other words, the validity of
\begin{align}\label{rombo}
\frac xy=z \quad \Otto \quad x=yz 
\end{align}

It is important to notice from the start that the equivalence above, read in set theoretical terms, says that the division operation has a value if and only if

\begin{align}\label{subset}
x \subseteq y.
\end{align}

 A major difference with respect to ordinary algebra is that the values $z$ are (almost) \emph{always} undetermined (not only for null $y$). In fact, $\frac xy$, thus $z$, is any set that intersected with $y$ results in $x$:
 
\begin{align}\label{division}
\frac xy \: y=x.
\end{align}

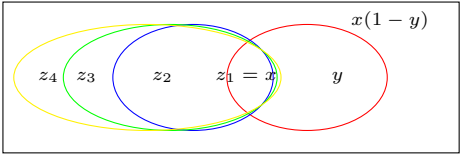
\begin{figure}[!bht]
\begin{tikzpicture}

\draw [color=red](05,01) ellipse (30pt and 20pt);    

\draw [color=blue] (03.5,01) ellipse (30pt and 20pt);    

\draw [color=green] (03.2,01) ellipse (40pt and 20pt);    

\draw [color=yellow] (02.9,01) ellipse (50pt and 20pt);    

\draw (01,00) rectangle (07,02);
\node at (06.09,02.0) [below] {\scriptsize{$x(1-y)$}};
\node at (05.4,01.2) [below] {\scriptsize{$y$}};
\node at (04.2,01.2) [below] {\scriptsize{$z_1=x$}};
\node at (03.1,01.2) [below] {\scriptsize{$z_2$}};

\node at (02.1,01.2) [below] {\scriptsize{$z_3$}};
\node at (01.6,01.2) [below] {\scriptsize{$z_4$}};

\end{tikzpicture}
\caption{$z$ is undetermined.}
\label{fig:division}
\end{figure}

 There are in general infinitely many possible values of $z$ which are suitable for the goal, more precisely, every set $z$ ranging from $x$ to $x+(1-y)$ will work. This reflects exactly Boole's words ``\emph{all}, \emph{some}, or \emph{none}'', concerning the members of $(1-x)(1-y)$ (p.90).\footnote{Of course, which one of the three cases will depend on the different individual examples. What we give here is a general theory designed by abstraction from single cases.} See the following Euler-Venn diagrams in {Figure} \ref{fig:division} for a visual clarification of this fact.\footnote{T. Hailperin, one of the very few authors, to our knowledge, suggesting a possible rigorous treatment of the division, shares with us the same starting point (\cite{H76b}. pp.70--77). But his approach is very different from ours: he develops a rigorous system with techniques of modern algebra that are external to Boole's conceptual background. We rather prefer to clarify the notion of division Boole worked with by shaping his fundamental intuitions through a suitable extension of the pseudo-binary calculus. In particular, for Hailpering $\frac00$ and $\frac10$ denote algebraic abstract entities, whereas they are  ordinary numbers in our approach.}


If we want to simulate the \emph{set-theoretical} calculus enriched with the division operator by a suitable  \emph{quantitative} algebraic calculus, we must therefore preserve somehow this idea of \emph{multi-valuedeness}. A possible way out is to say that any operation satisfying the three requirements below is a good interpretation of the division operation. In the following, $p$ and $q$ are arbitrary rational numbers and $div$ is used to denote the ordinary division between rational numbers (to avoid ambiguity with $/$):

\begin{enumerate}
\item $\frac pq=div(p,q)$, for $p\neq0\neq q$;
\item $\frac p0\in\IQ\setminus\{0,1\}$, whenever $p\neq0$;
\item $\frac 00\in\{0,1\}$.
\end{enumerate}

Notice that the second condition is in accordance with Boole's claim that $\frac10$ fails to fulfill the duality law, whereas the third condition meets his idea of $\frac00$ as an undetermined quantity that nevertheless should always represent a class.\footnote{Beside the conditions 1--3, the reader may reasonably ask for a forth natural requirement, that is, $\frac k0=k\frac10$ for $k\neq 0$. The reader is of course allowed to select a function respecting this further restriction, however this is not relevant for our discussion here.}

%
In this approach $\frac 00$ and  $\frac k0$  are admissible expressions denoting quantities, so the realm of fractions of standard algebra is extended.

Consider now the generalization of (\ref{rombo}):

\begin{align}\label{rombogen}
\frac st=u \quad \Otto \quad s=tu
\end{align}


\noindent
for {arbitrary} terms $s,t,u$. This holds in ordinary algebra only if $t\neq 0$, otherwise it has no meaning at all. On the contrary, in our new calculus this equation is valid also for $s=t=0$, whereas it is not valid for $s \neq 0$ and $t = 0$.  This last fact can be motivated in set theoretical terms by noticing that $\frac s 0$ has a meaning only when $s$ is a subset of the empty set, hence $\frac s 0$ is set theoretically not allowed when $s \neq 0$. On the other hand, $\frac 00$ is perfectly acceptable since  $\eps\sbsq\eps$ and it satisfies the fundamental equivalence (\ref{rombo}).


A direct consequence valid in the same domain (hence, in particular, for $s=0=t$) is the generalization of (\ref{division}),

\begin{align}\label{circledS}
\frac st \: t=s 
\end{align}

\noindent
(just substitute $u$ by $\frac st$ in virtue of $\frac st=u$): set theoretically, whenever $s$ is a subset of $t$, then $\frac st$ is a set that intersected with $t$ gives $s$.

The involved validity condition is then:
\begin{align}\label{triangle}
\textrm{\emph{no denumerator vanishes unless the corresponding numerator does the same}}
\end{align}
This condition extends the usual one required in ordinary algebra:
\begin{align}\label{triangle_triangle}
\textrm{\emph{no denumerator vanishes}}
\end{align}


It is then not difficult to see that under (\ref{triangle}), or in the case $u$ satisfies the law of duality (as it will be the case below), also the following law from ordinary algebra is preserved:

\begin{align}\label{odot}
\frac st \: u=\frac{su}{tu} \: u.
\end{align}

This relation will be soon of fundamental importance. 

However, here are two cases for which a larger domain of validity is not granted. The first one is particularly interesting:

\begin{align}\label{circledS_circledS}
\frac{st}t=s
\end{align}


\noindent
It fails to be true for $t=0$, $\frac00:=i$ and $s\neq i$.
This failure, however unpleasant at first sight, perfectly agrees with the fact that even if $st$ is always a subset of $t$, this does not imply that $s$ is the only set that intersected with $t$ gives us $st$.

Analogue observations apply to the equations
$$\frac st=\frac{su}{tu} \qquad \frac st \frac uv=\frac{su}{tv}$$
that do not necessarily hold: neither set-theoretically, nor numerically under (\ref{triangle}).

\vspace{0,2cm}

Nevertheless, the given definition of division suffices to prove our main requirement: the extension of Proposition \ref{developments-intersections-equivalence-a} to the division operator.

\begin{proposition}\label{developments-intersections-equivalence-b}
\emph{For functional terms $f(x,y)$ in the full algebraic language $\{+,\times,-,/,0,1\}$, the} method of the developments \emph{and the} method of the intersections \emph{are equivalent.}
\end{proposition}
First of all let us observe that the methods of the developments and the proof of Proposition \ref{developments} can be immediately transferred to the division. The proof of that proposition is in fact completely independent of the type of $f$, the only delicate point is to ensure that all possible fractional terms are always admissible, but this is indeed the case in our system.

Therefore, it remains just to complete the proof of Proposition
\ref{developments-intersections-equivalence-a} with the induction step:

\begin{itemize}
\item $f(x,y):=\frac{f_1(x,y)}{f_2(x,y)}$. Of course $f(i,j)=\frac{f_1(i,j)}{f_2(i,j)}$. By IH $f_l(x,y)c_{ij}$ and $f_l(i,j)c_{ij}$ coincide up to algebraic transformations for $l=1,2$. Then, by (\ref{odot}), since $c_{ij}$ satisfies the duality law, $f(x,y)c_{ij}:=\frac{f_1(x,y)}{f_2(x,y)}c_{ij}=\frac{f_1(x,y)c_{ij}}{f_2(x,y)c_{ij}}c_{ij}=\frac{f_1(i,j)c_{ij}}{f_2(i,j)c_{ij}}c_{ij}=\frac{f_1(i,j)}{f_2(i,j)}c_{ij}:=f(i,j)c_{ij}$.
\end{itemize}

The proof is now complete.

\vspace{0,5cm}

But this is not the end of the story. It is interesting to notice that by mere algebraic calculations an alternative justification of Boole's interpretation of division can be given. Such justification relies only on the syntactic rule (\ref{rombo}) and on an application of the method of the intersection to $\frac xy$:


 
\begin{align}\label{paragrafo_2}
\frac xy\, = \,\frac xy xy+\frac xy x(1-y)+\frac xy (1-x)y+\frac xy (1-x)(1-y).
\end{align}


Let now $\frac xy\, = \, z$ for some $z$. We obtain from (\ref{paragrafo_2})
$$\frac xy\, = \, z xy+z x(1-y)+z (1-x)y+z (1-x)(1-y) .$$
Since (\ref{rombo}) holds by assumption, we have
\\
$$\frac xy\, = \, z yzy+z yz(1-y)+z (1-yz)y+z (1-x)(1-y)\, = \, $$
\begin{align}\label{paragrafo2}
yz+yz(1-y)+z (1-yz)y+z(1-x)(1-y).
\end{align}
The variable $z$ is present in all the four addenda, but only in the last it plays an essential role, as the following calculations show:

\begin{itemize}
\item $yz=xy$, since from $\frac xy\, = \, z$ follows that $x=yz$, and so  $xy=y^2z=yz$;
\item $yz(1-y)=zy(1-y)=z\cdot 0=0$
\item $z(1-yz)y\, = \, (z-yz)y\, = \, z(1-y)y=0$
\end{itemize}

So far each addendum coincides with the corresponding one in the development (\ref{paragrafo}). 

$z$ remains only in the forth addendum $z(1-x)(1-y)$.
Since $z$ is an unknown parameter, the situation is exactly the one described by Boole's expression ``\emph{all}, \emph{some}, or \emph{none}'' referred to the members of $(1-x)(1-y)$ (p.90).

In conclusion, the use of the intersections yields the same intuitive interpretations of $\frac10$ and $\frac00$ suggested by Boole. And since $x$ and $y$ range over all possible classes, this proof has the highest degree of generality.

As a side remark, we observe that an hypothetical application of (\ref{circledS_circledS}), i.e. $\frac{st}t=s$, would on the contrary imply syntactically the unpleasant result $\frac xy=x$. First of all, observe that under  (\ref{division}) $x=xy$, since
$x=xx=\frac xy y x=  \frac xy y^2 x = \frac xy y yx = xyx = xy$. Hence, by (\ref{circledS_circledS}), $\frac xy=\frac{xy}y=x$. Although this result would deliver one of the suitable outputs, the multi-valudeness of the division would be lost in this way. But, as we have seen, $\frac{st}t=s$ is not valid under (\ref{triangle}) for every allowed interpretation of the division!

%
%

\vspace{1cm}

We conclude our analysis of the division by an inspection of its pseudo-binary table:

\vspace{0,2cm}

\begin{tabular}{p{0cm}p{0cm}p{0,5cm}p{2cm}}
$A$ & $B$ & & $A/B$\\
1 & 1 & & 1\\
1 & 0 & & $q\neq 0,1$\\
0 & 1 & & 0\\
0 & 0 & & 0 or 1\\
\end{tabular}

\vspace{0,2cm}

This table shows how the use of $\eps$ and $U$ reflects the behaviour of the division over all possible sets $A$ and $B$, analogously to what happens for the other operations! Always according to (\ref{rombo}) and (\ref{circledS}), the idea is that $A/B$ is defined only when $A\sbsq B$ and in this case $A/B$ is such that 
$$A=\frac AB \cap B.$$

As to the first line of the table above, $U=Z\cap U$ if and only if $Z=U$.
The second line instead originates no set, as $U$ is not a subset of $\eps$, and therefore it is set theoretically not allowed. 
As to the third line, the intersection of any set $Z$ with a non empty set results in $\eps$ if and only if $Z=\eps$.
For the fourth, we observe that $\eps=Z\cap\eps$ holds for all sets, in particular for $Z=\eps$ or $Z=U$. This is in perfect agreement with the fact that all sets $Z$ ranging from $A$ to $A+(1-B)$  are valid solutions. If only $\eps$ and $U$ are available, then  $Z=A=\eps$ or $Z=A+(1-B)=\eps+(U-\eps)=U$.
\vspace{0,2cm}

\section{Conclusion: the correspondence revisited}    

John Corcoran in his introduction to \emph{The Laws of Thought}, Prometheus Books 2003, says that  
\begin{quote}
Boole is one of the most misunderstood of the major philosophers of logic. He gets criticized for things he did not do, or did not do wrong. He never confused logic with psycology. He gets credits for things he did not do, or did not do right. He did not write the first book on pure mathematics and he did not devise `Boolean algebra'. Even where there is no question of blame or praise, his ideas are misdescribed or, worse, ignored.(pag xxix-xxx) 
\end{quote}
We do agree with Corcoran.  In the present paper we have tried to present the basis of Boole's logical theory as it emerges from \emph{The Laws of Thought} with its good and bad features. To this aim we have proved that the method of developments, once it is seen through the method of intersections, shows its clear set-theoretical meaning and it shows how natural and effective was Boole's idea of decomposing a concept in the way he did. In our analysis, it emerges that
the operations that Boole denotes by $+$ and $-$ are neither the standard set-theoretical union (between disjoint sets or not) nor the set-theoretical difference.
$+$ and $-$ are partial operations.\footnote{ Corcoran in his introduction agrees on this point, see pag xxx-xxxi).} What emerges from our work is that those operations are not defined, set-theoretically, exactly when the corresponding algebraic operations assume values outside $\{0, 1\}$.
 The division operation in the algebra of logic is again a partial operation, as it is in quantitative algebra.
\\
\\
 {Furthermore, we have tried to provide a mathematical interpretation justifying the logical conception of Boole's division and his intuitive treatment of unusual coefficients such as $\frac10$ and $\frac00$.}
\\
\\
At each stage of our presentation of Boole's theory we have insisted in checking if the alleged correspondence between algebra and logic were achieved or not, or in which degree. Do they share the same universal laws?
This point is of fundamental importance, since an essential part of Boole's machinery relies on the perfect formal correspondence between them so that we should nowadays rather speak of two different interpretations of a formally unique calculus.  As Boole himself says, when proving the validity of a logical law or of an argument through symbolic calculations, one would be allowed to \emph{suspend} the logical interpretation in all intermediate steps and rely {only} on the application of ordinary algebraic rules plus the law of duality. The {presumed} existence of a perfect formal correspondence between the logical and the quantitative pseudo-binary calculi would always guarantee the correctness of the proof:

\begin{quote}
It has been seen, that any system of propositions may be expressed by equations involving symbols $x$, $y$, $z$, which, whenever interpre- tation is possible, are subject to laws identical in form with the laws of a system of quantitative symbols, susceptible only of the values 0 and 1 [...]. But as the formal processes of reasoning depend only upon the laws of the symbols, and not upon the nature of their interpretation, we are permitted to treat the above symbols, $x$, $y$, $z$, as if they were quantitative symbols of the kind above described. \emph{We may in fact lay aside the logical interpretation of the symbols in the given equation; convert them into quantitative symbols, susceptible only of the values 0 and 1; perform upon them as such all the requisite processes of solution; and finally restore to them their logical interpretation.} (pp.69--70)
\end{quote}

Actually, despite Boole's enthusiatic slogan about the correspondence of the two calculi, he seems perfectly aware of the irreducible discrepancy produced by the need of interpretational restrictions for his logical calculus.
Nevertheless, instead of being worried about it, he thinks (quite unexepctedly, we would say) that he can take advantage of it:


\begin{quote}
The processes to which the symbols $x,y,z$ regarded as quantitative and of the species above described, are subject, are not limited by those conditions of thought to which they would, if performed upon purely logical symbols, be subject, and a freedom of operation is given to us in the use of them, without which, the inquiry after a general method in Logic would be a hopeless quest. (p.70)
\end{quote}

We see here somehow a \emph{sudden change of perspective}. Although Boole had insisted, as much as possible, on the perfect correspondence between the algebraic and the logical caluli, in the end his \emph{true} idea becomes that of \emph{substituting} the logical calculus by the quantative pseudo-binary calculus, taking advantage of the freedom of interpretation in the latter.  The subtle point to appreciate here is that such substitution is recommended because it is not vacuous:  the two calculi  do not coincide and one is more efficient than the other!

\vspace{0,2cm}

%

\end{document}